\documentclass[12pt]{amsart}
\usepackage[english]{babel}
\usepackage{amssymb}
\usepackage{amsfonts}
\usepackage{amsmath}
\usepackage{bm}
\usepackage{latexsym}
\usepackage{epsfig}
\numberwithin{equation}{section}
\usepackage{amssymb}

\newtheorem{theorem}{Theorem}[section]

\newtheorem{lemma}[theorem]{Lemma}
\newtheorem{proposition}[theorem]{Proposition}
\newtheorem{definition}[theorem]{Definition}

\newtheorem{remark}[theorem]{Remark}

\numberwithin{equation}{section}
\def\beq{\begin{equation}}
\def\eeq{\end{equation}}
\def\ben{\begin{enumerate}}
\def\een{\end{enumerate}}

\def\fa{{\mathfrak a}}
\def\fb{{\mathfrak b}}
\def\cF{{\mathcal F}}

\def\cH{{\mathcal H}}

\def\cM{{\mathcal M}}

\newcommand{\RR}{{\mathbb R}}
\newcommand{\CC}{{\mathbb C}}

\newcommand{\ZZ}{{\mathbb Z}}

\newcommand{\NN}{{\mathbb N}}

\newcommand{\ol}{\overline}

\usepackage{color}
\usepackage[normalem]{ulem}
\date{23 May 2018}

\begin{document}

\title[Generalized Fock spaces and the Stirling numbers]{Generalized Fock spaces and the Stirling numbers}

\author[D. Alpay]{Daniel Alpay}
\address{(DA) Schmid College of Science and Technology\\
Chapman University\\
One University Drive
Orange, California 92866\\
USA}
\email{alpay@chapman.edu}
\author[M. Porat]{Motke Porat}
\address{(MP) Department of Mathematics\\
Ben-Gurion University of the Negev\\
P.O. Box 653,
Beer-Sheva 84105\\ Israel}
\email{motpor@gmail.com}
%
%
\maketitle

\begin{abstract}
The Bargmann-Fock-Segal space plays an important role in mathematical physics, and has been extended
into a number of directions. In the present paper we imbed this space into a Gelfand triple. 
The spaces forming the Fr\'echet part (i.e. the space of test functions)
of the triple are characterized both in a geometric way and in terms of the adjoint of 
multiplication by the complex variable, using the Stirling numbers of the second kind. The dual of
the space of test functions has a topological algebra structure, of the kind introduced and studied by the first named 
author and G. Salomon.
\end{abstract}


\section{Introduction}
The reproducing kernel Hilbert space 
$\cF_1$ 
of entire functions with reproducing kernel 
$e^{z\overline{w}}$ 
is associated with the names of Bargmann, 
Segal and Fock, and will be called in this paper
the Fock space (more precisely, it is the symmetric 
Fock  space associated with
$\mathbb C$, see 
\cite{MR0157250}). 
It plays an important role in stochastic processes, 
mathematical physics and quantum mechanics, for recent work on the topic see e.g. 
\cite{DHK,Hall1}. 
The space $\cF_1$ 
is isometrically included in the Lebesgue space 
of the plane with weight 
$dA(z):=\frac{1}{\pi}e^{-|z|^2}dxdy$, 
and a key feature of 
$\cF_1$ 
is that the adjoint of the operator of 
multiplication by the complex variable is 
the operator of differentiation. It is of interest to 
look at various generalizations of 
$\cF_1$.
One approach consists in slightly 
modifying the weight function, see e.g. the works 
\cite{MR728694,rosenblum_hermite,siso} 
and another line is to change the 
kernel (that is the norms of the monomials) in an appropriate way, for 
instance replacing the 
exponential by the Mittag-Leffler 
function in the case of the grey noise theory, see e.g. 
\cite{MR1124240}. Then too the weight is changed, but not always in an explicit way.
Here we consider the family 
$(\cF_m)_{m=1}^\infty$ 
of reproducing kernel Hilbert spaces 
with reproducing kernel
\begin{align}
k_m(z,\omega)=
\sum_{n=0}^\infty 
\frac{z^n\ol{\omega}^n}{(n!)^m},\quad m=1,2,\ldots
\label{kmmm}
\end{align}
The space 
$\cF_m$ 
can then be easily described as the space 
of all Taylor series of the form 
$f(z)=\sum_{n=0}^\infty f_nz^n$ 
for which 
$$\sum_{n=0}^\infty 
|f_n|^2(n!)^m<\infty.$$
For 
$m=1$, 
the space is equal to the classical Fock space, and the case $m=2$ was defined and studied in \cite{daf1}.\\\\
The main results are as follows: The first is a geometric characterization of the spaces $\cF_m$ in terms of a weight; the second main result is the characterization of $\cF_m$
in terms of the adjoint operator of multiplication by $z$, associated to the Stirling numbers of second kind, 
see (\ref{eq:23May18a}). This generalizes the well known case of $m=1$, and opens the ground for future applications such as interpolation and sampling theorems 
in the setting of $\cF_m$; see for instance the papers \cite{MR2672228,MR3558232} for the case of $\mathcal F_1$.
The third main result is obtaining a structure of topological algebra for the inductive limit of the dual of the space $\cap_{m=1}^\infty\mathcal F_m$. 
This allows us to work locally in a Hilbert space rather than in the non-metrizable space $\cup_{m\in\mathbb N}\mathcal F_{2-m}$.\\

The outline of the paper is as follows. In Section \ref{2} we review some facts on the Mellin transform.
In Section \ref{3}, using the Mellin transform, we give a geometric characterization of the spaces $\mathcal F_m$ for $m\in\mathbb N$. A characterization 
of these spaces in terms of the adjoint of the operator of multiplication by $z$ and using the Stirling numbers of the second kind
 is given in Section \ref{4}. A related Bargmann
transform is defined in Section \ref{5}. In Section \ref{6} we define a Gelfand triple in which we imbed the Fock space.
We observe that the intersection $\bigcap_{m=1}^\infty\cF_m$ is a nuclear space and its dual is an algebra of the type introduced in \cite{MR3404695}.

\section{Preliminaries}
\setcounter{equation}{0}
\label{2}
Let $(a,b)$ a open interval of the real line, and let $f$ and $g$ such that both 
$f(x)x^{c-1}$ and $g(x)x^{c-1}$ are summable on $[0,\infty)$ for $c\in(a,b)$.
The Mellin transform of 
$f$, 
denoted by 
$\cM(f)$, 
is given by
\[
\cM(f)(c):=\int_0^\infty x^{c-1}f(x)dx,
\quad c\in(a,b).
\]
In particular, the Mellin transform of the function 
$f_1(x)=e^{-x}$ 
is the Gamma function:
$$\cM(f_1)(c)=
\int_0^\infty x^{c-1}e^{-x}dx=\Gamma(c),
\quad c>0.$$
The Mellin convolution of $f$ and $g$ is defined by  
$$(f*g)(x):=\int_0^\infty
f(\frac{x}{t})g(t)\frac{dt}{t}
=\int_0^\infty f(t)
g(\frac{x}{t})\frac{dt}{t},
\quad x>0.$$

An important relation between the 
Mellin transform and the Mellin convolution, 
see e.g. \cite[Theorem 3]{MR1468369}, is given by
\begin{align*}
\cM(f*g)(c)=(\cM(f)(c))(\cM(g)(c)),
\quad c\in(a,b).
\end{align*}
\section{Geometric description of 
$\cF_m$}
\setcounter{equation}{0}
\label{3}
Recall that the Fock space $\cF_1$ consists of those entire functions 
$f$ 
for which 
$$\iint_{\CC}|f(z)|^2 e^{-|z|^2}dA(z)<\infty,$$
and is the reproducing kernel Hilbert space 
with reproducing kernel 
$e^{z\overline{w}}$.
In this section we give for $m=2,\ldots$ a geometric characterization for the space 
$$\cF_m=\left\{
f(z)=\sum_{n=0}^\infty a_n z^n
\text{ is entire with }
\sum_{n=0}^\infty |a_n|^2
(n!)^m<\infty\right\}$$
which is the reproducing kernel Hilbert space with reproducing kernel \eqref{kmmm}, when equipped with the inner product
\begin{align*}
\langle f,g\rangle_{\cF_m}:=
\sum_{n=0}^\infty f_n\ol{g_n}(n!)^m,
\,\text{where }
f(z)=\sum_{n=0}^\infty f_nz^n, \,\,
g(z)=\sum_{n=0}^\infty g_nz^n,
\end{align*}
for every 
$f,g\in\cF_m$. 
First, we use the properties of the 
Mellin transform to build 
the kernels 
$K_m(z)$, 
which are generalizations of the 
modified Bessel function of the second order,
also called the Macdonald function. 
Let $K_1(x)=e^{-x}$ and for every integer $m>1$ define the function
\begin{align}
\label{eq:22May18a}
K_m(x):=(K_1*\cdots*K_1)(x),
\quad x\in\RR_+
\end{align}
that is the function
$K_1(x)$ 
Mellin-convoluted 
$m$ 
many times with itself.
\begin{lemma} 
Let 
$m$ 
be an integer. The following properties hold: 
\begin{itemize}
\item[(1)] 
For $m>1$, 
the kernel 
$K_m$ 
has the integral representations
\begin{align}
\label{eq:7Apr18a}
K_m(x)=\int_0^\infty\cdots\int_0^\infty
\frac{e^{-\sum_{i=1}^{m-1}
x_i-\frac{x}{\prod_{i=1}^{m-1} x_i}}}
{\prod_{i=1}^{m-1}x_i}
dx_1\cdots dx_{m-1}
\end{align}
and
\begin{align}
\label{eq:7Apr18b}
K_m(x)=\int_{\RR}\cdots\int_{\RR}
e^{-\sqrt[m]{x}(\sum_{i=1}^{m-1}e^{t_i}
+e^{-\sum_{i=1}^{m-1}t_i})}dt_1\cdots dt_{m-1}.
\end{align} 
\item[(2)]
The function 
$K_m$ 
is monotone decreasing in 
$(0,\infty)$.
\item[(3)] 
The Mellin transform of 
$K_m$ 
is given by 
\[
\cM(K_m)(x)=\Gamma(x)^m,\quad x>0,
\]
and so 
\begin{align}
\label{eq:9Apr18a}
\int_0^\infty x^nK_m(x)dx
=(n!)^m,\quad n\in\NN.
\end{align}
\end{itemize}
\end{lemma}
\textbf{Proof.}
Part 1 is proved by induction on 
$m$: 
if 
$m=2$, 
we get
$$K_2(x)=\int_0^\infty 
e^{-x/t}e^{-t}
\frac{dt}{t}=\int_0^\infty 
\frac{e^{-x_1-
\frac{x}{x_1}}}{x_1}dx_1.$$
Suppose formula 
(\ref{eq:7Apr18a}) 
holds for 
$m$. Then
\begin{align*}
K_{m+1}(x)&=(K_m*e^{-t})(x)
=\int_0^\infty
K_m\left(\frac{x}{x_m}
\right)e^{-x_m}
\frac{dx_m}{x_m}
\\&=\int_0^\infty\cdots\int_0^\infty
\frac{e^{-\sum_{i=1}^{m-1}
x_i-\frac{\frac{x}{x_m}}
{\prod_{i=1}^{m-1} x_i}}}
{\prod_{i=1}^{m-1}x_i}
\frac{e^{-x_m}}{x_m}
dx_1\cdots dx_m 
\\&=\int_0^\infty\cdots\int_0^\infty 
\frac{e^{-\sum_{i=1}^m x_i-
\frac{x}{\prod_{i=1}^m x_i}}}
{\prod_{i=1}^m x_i} 
dx_1\cdots dx_m,
\end{align*}
i.e., 
(\ref{eq:7Apr18a}) 
holds for 
$m+1$
and hence for every 
$m>1$. 
Next, we use 
(\ref{eq:7Apr18a}) 
and the change of variables 
$s_i=\ln(x_i)$, $1\le i\le m-1,$
to obtain
\begin{align*}
K_m(x)=\int_{\RR}\cdots\int_{\RR}
e^{-\sum_{i=1}^{m-1}e^{s_i}-
\frac{x}{e^{\sum_{i=1}^{m-1}s_i}}}
ds_1\cdots ds_{m-1}
\end{align*} 
and by another change of variables 
$t_i=s_i-\ln(\sqrt[m]{x}), 
1\le i\le m-1,$ we get
$$K_m(x)=\int_{\RR}\cdots
\int_{\RR} e^{-\sqrt[m]{x}
(\sum_{i=1}^{m-1}e^{t_i}+
e^{-\sum_{i=1}^{m-1}t_i})}
dt_1\cdots dt_{m-1}.$$
From the representation 
(\ref{eq:7Apr18a}) 
it is easily seen that 
$K_m(x)$ 
is a monotone decreasing function. Finally, the Mellin transform of 
$K_m$ 
is given by
\begin{align*}
\cM(K_m)(c)=
\cM(f_1)(c)\cdots\cM(f_1)(c)
=(\Gamma(c))^m,\quad c>0,
\end{align*}
therefore 
$$\int_0^\infty x^{c-1}K_m(x)dx
=(\Gamma(c))^m,\quad c>0.$$
For $c=n+1$, we have
$$\int_0^\infty x^nK_m(x)dx
=(\Gamma(n+1))^m=(n!)^m.
\quad\blacksquare$$
In the special case $m=2$, we get that 
\[
K_2(x)=\int_\RR e^{-\sqrt{x}2\cosh(t)}dt,\quad x\in\RR_+
\]
is the Bessel function of the second kind, see 
\cite{daf1}. For an arbitrary $m>2$, the kernel $K_m(x)$ can be expressed in terms of the Meijer $G$-functions; 
\cite[Chapter 5]{MR0058756} for the latter. We now show how the generalized Fock spaces $\cF_m$ are obtained from the kernels $K_m(x)$ in a natural way. 
\begin{theorem}
For any integer 
$m\ge1$, 
the space 
$\cF_m$ 
is equal to the space of all entire functions 
$f:\CC\rightarrow\CC$ 
satisfying the condition
\begin{align}
\label{eq:7Apr18c}
\iint_{\CC}|f(z)|^2
K_m(|z|^2)dA(z)<\infty.
\end{align}
Moreover, the inner product of $\cF_m$ is given by
\begin{align*}
\frac{1}{\pi}\iint_{\CC}
f(z)\ol{g(z)}K_m(|z|^2)dA(z)
=\sum_{n=0}^\infty 
f_n\ol{g_n}(n!)^m,\quad f,g\in\cF_m,
\end{align*}
and 
$\cF_m$ 
has the orthonormal basis 
$\left\{\frac{z^n}{(n!)^{m/2}}
\right\}_{n=0}^\infty$.
\end{theorem}
\textbf{Proof.}
A straightforward computation shows that
\begin{align*}
\iint_{\CC}
z^n\ol{z}^kK_m(|z|^2)dA(z)
&=\int_0^\infty 
\int_0^{2\pi}r^ne^{in\theta}
r^ke^{-ik\theta}
K_m(r^2)rd\theta dr
\\&=\int_0^{2\pi}
e^{i(n-k)\theta}d\theta
\int_0^\infty r^{n+k+1}K_m(r^2)dr
\\&=2\pi\delta_{n,k}
\int_0^\infty r^{2n+1}K_m(r^2)dr
\\&=2\pi\delta_{n,k}
\int_0^\infty u^nK_m(u)
\frac{du}{2}\\
&=\pi (n!)^m\delta_{n,k}.
\end{align*}
Let 
$f=\sum_{n=0}^\infty 
f_nz^n$ 
and 
$g=\sum_{n=0}^\infty g_nz^n$ 
be entire functions. Then
\begin{align*}
\pi \iint_{\CC}f(z)\ol{g(z)}
K_m(|z|^2)dA(z)&=
\sum_{n,k=0}^\infty 
f_n\ol{g_k}\iint_{\CC}z^n
\ol{z}^kK_m(|z|^2)dA(z)
\\&=\sum_{n,k=0}^\infty 
f_n\ol{g_k}\delta_{n,k}(n!)^m
=\sum_{n=0}^\infty
f_n\ol{g_n}(n!)^m,
\end{align*}
which implies that $f\in\cF_m$ if and only if condition  (\ref{eq:7Apr18c}) holds, i.e., 
\[
\frac{1}{\pi}
\sum_{n=0}^\infty|f_n|^2(n!)^m
=\iint_{\CC}
|f(z)|^2K_m(|z|^2)dA(z)<\infty
\] 
as wanted. Furthermore, the inner product in $\cF_m$ is then given by
$$\langle f,g\rangle_{\cF_m}=
\sum_{n=0}^\infty f_n\ol{g_n}(n!)^m
=\frac{1}{\pi}\iint_{\CC}f(z)
\ol{g(z)}K_m(|z|^2)dA(z). 
\quad\blacksquare$$
In the case $m=2$,  similar yet different spaces related to other families of orthogonal polynomials, appear in  
\cite[Lemma 4]{Karp1} and \cite{Karp2}, 
\begin{remark}
Let 
$0<\epsilon<1$.
Then 
$\frac{\epsilon}{n!}<1$ 
for every 
$n\ge0$ 
and hence
\begin{align*}
\sum_{m=1}^\infty 
\epsilon^mk_m(z,\omega)
&=\sum_{m=1}^\infty \epsilon^m
\left(\sum_{n=0}^\infty 
\frac{z^n\ol{\omega}^n}{(n!)^m}\right)
=\sum_{n=0}^\infty \left(
\sum_{m=1}^\infty\left(
\frac{\epsilon}{n!} 
\right)^m\right)z^n\ol{\omega}^n
\\&=\sum_{n=0}^\infty 
\frac{\epsilon}{n!}
\left(
\frac{1}{1-\frac{\epsilon}{n!}}
\right) 
z^n\ol{\omega}^n
=\epsilon\cdot\sum_{n=0}^\infty 
\frac{z^n\ol{\omega}^n}{n!-\epsilon}
\end{align*}
and
\begin{align*}
\sum_{m=1}^\infty 
\frac{\epsilon^m}{m!}
k_m(z,\omega)=\sum_{n=0}^\infty
\left( e^{\frac{\epsilon}{n!}}
-1\right)z^n\ol{\omega}^n.
\end{align*}
\end{remark}
\section{Operator theoretic description of $\cF_m$}
\setcounter{equation}{0}
\label{4}
Denote by $\fa$ the operator of  multiplication by $z$ by $\fb$ 
differentiation by $z$, i.e., 
$\fa=M_z$ 
and 
$\fb=\frac{\partial}{\partial z}$. 
Both $\fa$ 
and
$\fb$ 
are defined on polynomials and more generally on 
entire functions. They satisfy the familiar commutation relation
$$[\fb,\fa]=\fb\fa-\fa\fb=I.$$
In the Fock space 
$\cF_1$, $\fa$ and $\fb$ are unbounded operators, and satisfy
$$\fa^*=\fb
\quad\text{ and }\quad 
\fb^*=\fa.$$
This relation is very important, 
as the Fock space is the only space of 
entire functions for which 
$\fa$ 
and 
$\fb$ are adjoint to each other, see 
\cite{MR0157250}.
We generalize this result by presenting  
a relation between the operators 
$\fa$ 
and 
$\fb$ 
in the space 
$\cF_m$. 
That gives us another characterization of the space $\cF_m$.
We first introduce the Stirling numbers of the second kind $S(k,n)$, which appear naturally in the theory of ordering bosons. 
\begin{definition}
[Stirling numbers of the second kind]
For $k\in\NN_0$ and $n\in\NN_0$, the numbers $S(k,n)$ are defined by the recurrence formula 
$$S(k,n)=
nS(k-1,n)+S(k-1,n-1),\quad k,n\ge1$$
with the initial values 
$S(k,0)=\delta_{k,0}$
and 
$S(k,n)=0$ if $k<n$.
\end{definition}
It is well known, see \cite{MR2862989,MR2479303}, that  $$(\fa\fb)^k=\sum_{n=1}^{k} S(k,n)\fa^{n}\fb^{n},
\quad k\ge 1$$
and this operator is called the Mellin derivative operator of order $k$ (with $c=0$), 
see \cite[Lemma 9]{MR1468369}.
\begin{theorem}
\label{thm:11Apr18a}
Let $m\ge1$ be an integer. The operators $\fa$ and $(\fb\fa)^{m-1}\fb$ are closed densely defined operators on the space $\cF_m$ 
and their domains coincide
$$Dom(\fa)=Dom((\fb\fa)^{m-1}\fb)=D,$$ 
where
\begin{align}
\label{eq:10Apr18a}
D=\left\{f(z)=\sum_{n=0}^\infty f_nz^n:
\sum_{n=0}^\infty|f_n|^2(n!)^mn^m<\infty\right\}
\subseteq\cF_m.
\end{align}
Moreover, the adjoint operator of
$\fa$ in $\cF_m$
is given by
$$\fa^*=(\fb\fa)^{m-1}\fb,\quad\text{with}\quad
Dom(\fa^*)=Dom\left((\fb\fa)^{m-1}\fb\right)=D.$$
Furthermore, let 
$\cH$ 
be a Hilbert space of entire functions in which the polynomials are dense, and let $m\in\NN$.
 If the adjoint operator of 
$\fa$ 
in $\cH$
is equal to the operator 
$(\fb\fa)^{m-1}\fb$, 
i.e., if
\begin{align}
\label{eq:23May18a}
(M_z)^*=
\frac{\partial}{\partial z}
\left[\sum_{n=1}^{m-1}S(m-1,n)z^n
\frac{\partial^n}{\partial z^n}\right],
\end{align}  
then 
$\cH=\cF_m$ and there exists 
$c>0$ 
for which 
$$\langle f,g\rangle_{\cH}=
c\cdot\langle f,g\rangle_{\cF_m},
\quad \forall f,g\in\cH.$$ 
\end{theorem}
\textbf{Proof.}
It is easy to see that 
$\fa$ 
and 
$(\fb\fa)^{m-1}\fb$ are closed, densely defined operators on $\cF_m$. If $f(z)=\sum_{n=0}^\infty
f_nz^n\in\cF_m$, then 
\begin{align*}
f\in Dom(\fa)&\iff 
\fa f=\sum_{n=0}^\infty f_nz^{n+1}\in\cF_m
\iff \sum_{n=0}^\infty 
|f_n|^2((n+1)!)^m<\infty
\end{align*}
and
\begin{align*}
f\in Dom((\fb\fa)^{m-1}\fb)&\iff
(\fb\fa)^{m-1}\fb f=\sum_{n=1}^\infty 
f_nn^mz^{n-1}\in\cF_m
\\&\iff \sum_{n=1}^\infty 
|f_n|^2n^{2m}((n-1)!)^m=
\sum_{n=1}^\infty |f_n|^2(n!)^mn^m<\infty.
\end{align*}
Therefore, 
$Dom(\fa)=Dom((\fb\fa)^{m-1}\fb)=D$
as in (\ref{eq:10Apr18a}).
Next, if 
$$g(z)=\sum_{n=0}^\infty 
g_nz^n\in Dom(\fa^*),$$  
there exists 
$$h(z)=\sum_{n=0}^\infty h_nz^n\in\cF_m$$
such that 
$\langle \fa f,g\rangle_{\cF_m}
=\langle f,h\rangle_{\cF_m}$
for every $f\in Dom(\fa)$. 
In particular, for 
$f(z)=z^n\, (n\ge0)$, 
we get
\begin{align*}
\ol{g_{n+1}}((n+1)!)^m=
\langle z^{n+1},g\rangle_{\cF_m}
=\langle z^n,h\rangle_{\cF_m}=\ol{h_n}(n!)^m
\end{align*}
and hence $h_n=g_{n+1}(n+1)^m$ for every $n\ge0$. Thus,
\begin{align*}
h\in\cF_m&\Longrightarrow \sum_{n=0}^\infty |h_n|^2(n!)^m=\sum_{n=0}^\infty |g_{n+1}|^2(n+1)^{2m}(n!)^m<\infty
\\&\Longrightarrow \sum_{n=1}^\infty |g_n|^2 (n!)^mn^m<\infty
\Longrightarrow g\in D,
\end{align*}
hence $Dom(\fa^*)\subseteq D$. Finally, if 
$g\in D=Dom((\fb\fa)^{m-1}\fb)$, 
then 
\begin{align*}
\langle f,(\fb\fa)^{m-1}
\fb g\rangle&=
\left\langle \sum_{n=0}^\infty 
f_nz^n,
\sum_{n=0}^\infty
(n+1)^mg_{n+1}z^n\right\rangle
=\sum_{n=0}^\infty 
f_n(n+1)^m\ol{g_{n+1}}(n!)^m
\\&=\sum_{n=0}^\infty 
f_n\ol{g_{n+1}}((n+1)!)^m=
\left\langle \sum_{n=0}^\infty 
f_nz^{n+1},\sum_{n=0}^\infty 
g_nz^n\right\rangle=\langle
\fa f,g\rangle,
\end{align*}
for every 
$f\in D=Dom(\fa),$
which proves that
$g\in Dom(\fa^*)$. Therefore, $D\subseteq Dom(\fa^*)$ and hence 
$Dom(\fa^*)=D$.
By the previous calculation we also know that  
$\fa^*=(\fb\fa)^{m-1}\fb$. Now  suppose that
$\cH$ 
is a Hilbert space which contains all polynomials, and such that 
$$\fa^*=(\fb\fa)^{m-1}\fb$$ in $\cH$. 
Then for every 
$f\in Dom(\fa)\cap\cH$
and
$g\in Dom((\fb\fa)^{m-1}\fb)\cap\cH$, 
\begin{align}
\label{eq:10Apr18b}
\langle \fa f,g
\rangle_{\cH}=\langle f,(\fb\fa)^{m-1}
\fb g\rangle_{\cH}
\end{align} 
and as both 
$Dom(\fa)$ and
$Dom((\fb\fa)^{m-1}\fb)$
contain all polynomials, we apply 
(\ref{eq:10Apr18b}) for the choice
$f(z)=z^l,g(z)=z^k\, (k,l\ge0)$, thus
\begin{align*}
\langle z^{l+1},z^k\rangle_{\cH}&=
\langle \fa f,g\rangle_{\cH}=\langle 
f,(\fb\fa)^{m-1}\fb g
\rangle_{\cH}
\\&=\langle 
z^l,k^mz^{k-1}\rangle_{\cH}=
k^m\langle z^l,z^{k-1}
\rangle_{\cH},\quad k,l\ge0.
\end{align*}
We now prove  by induction that for every 
$k\ge0$ and $l\ge k$, 
$$\langle z^{l+1},z^k\rangle_{\cH}=0:$$ 
\begin{itemize}
\item
If 
$k=0$, 
we know that 
$\langle z^{l+1},1\rangle_{\cH}=0$ 
for every 
$l\ge0$.
\item
Assume that for some 
$k\ge0$ 
we have 
$\langle z^{l+1},z^k\rangle_{\cH}=0$ 
for every 
$l\ge k$. 
Therefore, 
$\langle z^{l+2},z^{k+1}\rangle_{\cH}=(k+1)^m
\langle z^{l+1},z^k\rangle_{\cH}=0$
for every 
$l\ge k$, 
which means that 
$$\langle z^{l+1},z^{k+1}\rangle_{\cH}=0$$ 
for every 
$l\ge k+1$, as wanted.
\end{itemize}
Thus the family 
$\{z^k\}_{k=0}^\infty$ is orthogonal in $\cH$ and one can easily see  that 
\[
\langle z^k,z^k\rangle_{\cH}
=k^m\langle z^{k-1},z^{k-1}\rangle_{\cH},\quad \forall k\ge 1,
\]
which implies that
$$\langle z^k,z^k\rangle_{\cH}
=(k!)^m\langle 1,1\rangle_{\cH}.$$
To conclude, if 
$f(z)=\sum_{k=0}^\infty f_kz^k$ and $g(z)=\sum_{k=0}^\infty g_kz^k\in\cH$, 
then
$$\langle f,g\rangle_{\cH}=
\sum_{k,l=0}^\infty f_k\ol{g_l}
\langle z^k,z^l\rangle_{\cH}
=\sum_{k=0}^\infty f_k\ol{g_k}
(k!)^m\langle1,1\rangle_{\cH},$$
i.e., the inner product in 
$\cH$ 
is equal to the one in 
$\cF_m$, 
up to a positive multiplicative constant 
$c=\langle1,1\rangle_{\cH}$. As 
$\cH$ 
is a Hilbert space which contains all 
the polynomials, it follows that 
$$\cH=\left\{ f=\sum_{n=0}^\infty 
f_nz^n:\langle f,f\rangle_{\cH}=
c\sum_{n=0}^\infty
|f_n|^2(n!)^m<\infty\right\}=\cF_m.
\quad\blacksquare$$

In the previous theorem 
we proved that $\cF_m$ 
is the only  Hilbert space which contains all polynomials and in which the adjoint operator of 
$\fa=M_z$ is equal to the operator 
$$(\fb\fa)^{m-1}\fb=\frac{\partial}{\partial z}
\left[\sum_{n=1}^{m-1}S(m-1,n)z^n
\frac{\partial^n}{\partial z^n}\right].$$
One can see that we have the relations 
$$\fb^n\fa=\fa\fb^n+n\fb^{n-1}
\quad\text{and}\quad \fb\fa^n=\fa^n\fb+n\fa^{n-1}$$ 
for every $n\in\NN$, 
and in particular the operators $\fa$ and $\fa^*$ do not satisfy the commutation relation. However we have the following result.
\begin{proposition}
The commutator of $\fa$ and $\fa^*=(\fb\fa)^{m-1}\fb$ is equal to
\begin{align}
\label{eq:11Apr18b}
[\fa^*,\fa]=
I+\sum_{n=1}^{m-1}(n+1)S(m,n+1)\fa^n\fb^n.
\end{align}
\end{proposition}
\textbf{Proof.}
As
$$\fa^*=(\fb\fa)^{m-1}\fb=
\fb\sum_{n=1}^{m-1} S(m-1,n)\fa^{n}\fb^{n},$$
we have
\begin{align*}
[\fa^*,\fa]&=
\fb
\sum_{n=1}^{m-1} 
S(m-1,n)\fa^n\fb^n\fa-
\fa\fb\sum_{n=1}^{m-1} 
S(m-1,n)\fa^n\fb^n
\\&=\fb
\sum_{n=1}^{m-1}S(m-1,n)
\fa^n(\fa\fb^n+n\fb^{n-1})
-\fa\fb\sum_{n=1}^{m-1}
S(m-1,n)\fa^n\fb^n
\\&=(\fb\fa-\fa\fb)
\sum_{n=1}^{m-1}S(m-1,n)
\fa^n\fb^n+\fb\sum_{n=1}^{m-1}
nS(m-1,n)\fa^n\fb^{n-1}
\\&=\sum_{n=1}^{m-1}S(m-1,n)
\fa^n\fb^n+\sum_{n=1}^{m-1}nS(m-1,n)
(\fa^n\fb+n\fa^{n-1})\fb^{n-1}
\\&=\sum_{n=1}^{m-1}(n+1)
S(m-1,n)\fa^n\fb^n+
\sum_{n=1}^{m-1}n^2S(m-1,n)
\fa^{n-1}\fb^{n-1}
\end{align*}
and as $S(m-1,1)=S(m-1,m-1)=S(m,m)=1$, we have
\begin{align*}
[\fa^*,\fa]&=I+m\fa^{m-1}\fb^{m-1}
+\sum_{n=1}^{m-2}(n+1)
[S(m-1,n)+(n+1)S(m-1,n+1)]\fa^n\fb^n
\\&=I+m\fa^{m-1}\fb^{m-1}+
\sum_{n=1}^{m-2}(n+1)S(m,n+1)\fa^n\fb^n
\\&=I+\sum_{n=1}^{m-1}(n+1)
S(m,n+1)\fa^n\fb^n.\quad\blacksquare
\end{align*}
Sequentially, a straightforward calculation shows that
\begin{align*}
\|\fa f\|^2_{\cF_m}
=\|\fa^* f\|^2_{\cF_m}
+\|f\|^2_{\cF_m}+\sum_{k=1}^{m-1}\binom{m}{k}
\left[\sum_{n=0}^\infty|f_n|^2(n!)^mn^k\right]
\end{align*}
for every $f\in D$, which guarantees that all the terms in the identity are finite. It is tempting to write the last identity (with some abuse of notation) as
\begin{align*}
\|\fa f\|^2_{\cF_m}
=\|\fa^* f\|^2_{\cF_m}
+\|f\|^2_{\cF_m}+
\sum_{k=1}^{m-1}\binom{m}{k}
\langle f,(\fa\fb)^kf\rangle_{\cF_m},
\end{align*}
however 
$f\in D$ 
does not necessarily imply that 
$f\in Dom((\fa\fb)^k)$. 

Finally, we have the following relation between the operators $\fa,\fb$ and the family of spaces $(\cF_m)_{m\in\ZZ}$. 
For every $n\ge1$,
\begin{itemize}
\item 
the Fock space $\cF_1$ satisfies
\begin{align*}
\fa^n(\cF_1)\subseteq\cF_0
\quad\,\,\,\,\,\,\text{and}\quad \,\,\,\fb^n(\cF_1)\subseteq\cF_0,
\end{align*}
\item
if $m>1$, then
\begin{align*}
\fa^n(\cF_m)\subseteq\cF_{m-1}
\quad\text{and}\quad \,\,\,\fb^n(\cF_m)\subseteq\cF_m, 
\end{align*}
\item
if $m<1$, then
\begin{align*}
\fa^n(\cF_m)\subseteq\cF_m
\quad\quad\,\,\,\,\,\,\text{and}\quad
\,\,\,\fb^n(\cF_m)\subseteq\cF_{m-1}.
\end{align*}
\end{itemize} 
\begin{remark}
Unlike the situation in the the 
Fock space, where the adjoint of 
$\fb$ 
is equal to 
$\fa$, 
in the space 
$\cF_m$
the adjoint operator of 
$\fb$ 
is equal to
$$\fb^*\left(\sum_{k=0}^\infty f_kz^k\right)
=\sum_{k=0}^\infty 
\frac{f_k}{(k+1)^{m-1}}z^{k+1},$$
thus $\fb^*\ne\fa$ if $m>1$.

\end{remark}
\section{Generalized Bargmann Transform}
\setcounter{equation}{0}
\label{5}
 Recall that the normalized Hermite functions are defined by
$$\eta_n(t)=\frac{1}{\pi^{1/4}2^{n/2}
\sqrt{n!}}
e^{\frac{t^2}{2}}\left( e^{-t^2}
\right)^{(n)},
\quad n\in\NN_0.$$
The family $\{\eta_n\}_{n=0}^\infty$
is an orthonormal basis of the Lebesgue space ${\mathbf L}_2(\RR,dt).$
Furthermore, see \cite[p. 436]{MR1502747}, the $\eta_n$ are uniformly bounded
by some constant, i.e.,   
\begin{align*}
\exists C>0
\text{ such that 
$|\eta_n(t)|\le C$, 
for every 
$n\in\NN$ 
and 
$t\in\RR$.}
\end{align*}
Similarly to the symmetric Fock space associated with 
$\CC$, see e.g. 
\cite{MR0157250},
 that is 
$\cF_1$,
there is a a fourth characterization 
of the space 
$\cF_m$,
given by a mapping 
from 
${\mathbf L}_2(\RR,dt)$
into
$\cF_m$, presented in the following proposition.
\begin{proposition}
Let 
$m\ge2$. 
For every 
$t\in\RR$ and $z\in\CC$ define the function
\begin{align}
\label{eq:2Dec15a}
h_m(z,t):=\sum_{n=0}^\infty 
\frac{z^n}{(n!)^{m/2}}
\eta_n(t).
\end{align}
Then,
\begin{itemize}
\item[1.]
for every $t\in\RR$, the function $h_m(\cdot,t)$ 
is entire.
\item[2.]
 
$f\in \cF_m$ 
if and only if there exists 
$g\in {\mathbf L}_2(\RR,dt)$, 
such that
\begin{align}
f(z)=\int_{\RR} h_m(z,t)
g(t)dt=\langle g,
\overline{h_m(z,\cdot)}\rangle_{{\mathbf L}_2(\RR,dt)}.
\end{align}
\end{itemize}
\end{proposition}
\textbf{Proof.} 
Since the functions 
$\eta_n(t)$ 
are all bounded by 
$C$,
the sum in 
(\ref{eq:2Dec15a}) 
converges and so 
$h_m(\cdot,t)$ 
is entire. Next,
let
$f(z)=\langle g,
\overline{h_m(z,\cdot)}\rangle_{{\mathbf L}_2(\RR,dt)}$ for some $g\in {\mathbf L}_2(\RR,dt)$. Then, 
\begin{align*}
f(z)=\int_{\RR} \left(\sum_{n=0}^\infty      
\frac{z^n}{(n!)^{m/2}}
\eta_n(t)g(t)\right)dt=
\sum_{n=0}^\infty 
\frac{z^n}{(n!)^{m/2}}
\int_{\RR}\eta_n(t)g(t)dt.
\end{align*} 
As the system
$\{\eta_n\}_{n=0}^\infty$ forms an orthonormal basis of
${\mathbf L}_2(\RR,dt)$, 
we have Parseval's equality
$$\sum_{n=0}^\infty
\left| \int_{\RR}\eta_n(t)
g(t)dt\right|^2=
\int_\RR|g(t)|^2dt$$
and hence 
$f\in \cF_m$, 
since
$$\sum_{n=0}^\infty \left|
\frac{1}{(n!)^{m/2}}
\int_{\RR}\eta_n(t)g(t)dt
\right|^2(n!)^m=
\|g\|_{{\mathbf L}_2(\RR,dt)}^2<\infty.$$
Finally, let
$f\in \cF_m$. It can be written as
$f(z)=\sum_{n=0}^\infty a_nz^n$ 
with 
$\sum_{n=0}^\infty 
|a_n|^2(n!)^m<\infty$.
Setting
$$g(t)=\sum_{n=0}^\infty 
(n!)^{m/2}a_n\eta_n(t),$$ 
we observe that 
$$\|g\|_{{\mathbf L}_2(\RR,dt)}^2=
\sum_{n=0}^\infty 
|a_n|^2(n!)^m<\infty$$
and finally that
$$\langle h_m(z,\cdot),g\rangle_{{\mathbf L}_2(\RR,dt)}
=\sum_{n=0}^\infty 
\frac{z^n}{(n!)^{m/2}}(n!)^{m/2}a_n=f(z).
\quad\blacksquare$$
This characterization of 
$\cF_m$ 
motivates us to consider an 
associated Bargmann transform.
For any 
$g\in {\mathbf L}_2(\RR,dt)$
we define the Bargmann transform of $g$ to be
\begin{align*}
\mathfrak{B_m}(g):=
\sum_{n=0}^\infty 
\frac{z^n}{(n!)^{m/2}}
\int_{\RR}\eta_n(t)g(t)dt=
\langle g, \overline{h_m(z,\cdot)}\rangle_{{\mathbf L}_2(\RR,dt)}.
\end{align*} 
The mapping 
$\mathfrak{B_m}:{\mathbf L}_2(\RR,dt)\rightarrow\cF_m$ 
is unitary; it 
satisfies 
$$\mathfrak{B_m}(\eta_n)(z)=
\frac{z^n}{(n!)^{m/2}}
\quad\text{and}\quad\|g\|_{{\mathbf L}_2(\RR,dt)}
=\|\mathfrak{B_m}(g)\|_{\cF_m}$$
for every 
$g\in {\mathbf L}_2(\RR,dt)$.
\begin{remark}
In case where 
$m=1$, 
$\mathfrak{B_1}$ 
is the well known Bargmann transform 
and the function 
$h_1(z,t)$ can be written in closed form as
$$h_1(z,t)=e^{2tz-t^2-z^2/2}.$$
When $m>1$, finding an explicit closed formula for the function $h_m(z,t)$ might involve new generalizations of the exponential function. 
\end{remark}

\section{A Gelfand Triple associated to the family
$(\cF_m)_{m\in\mathbb Z}$}
\setcounter{equation}{0}
\label{6}
The reproducing kernel Hilbert spaces 
$\{\cF_m\}_{m=1}^\infty$,
starting from the Fock space 
$\cF_1$, 
form a decreasing sequence, i.e., 
$$\cF_1\supset \cF_2\supset...
\supset \cF_{m}\supset \cF_{m+1}\supset....$$
so it makes sense, in the spirit of the theory of Gelfand triples (as developed for instance in the books \cite{MR0435834,GS2_english}) to consider the intersection space
\[
\begin{split}
\cF&=\bigcap_{m=1}^\infty \cF_m\\
&=
\left\{ f=\sum_{n=0}^\infty a_n z^n 
\text{ such that }
\|f\|_m=\sum_{n=0}^\infty 
|a_n|^2(n!)^m<\infty,
\forall m\in\NN\right\}
\end{split}
\]
which consists of entire functions, and its dual. We consider the dual space of each 
$\cF_m$, with respect to the Fock space $\cF_1$.

\begin{lemma}
\label{lem:5Mar16a}
For every 
$m\ge1$, the dual space of 
$\cF_m$, 
with respect to 
$\cF_1$ is
$$\cF_{2-m}:=(\cF_m)^\prime
=\left\{ b=(b_n)_{n\in\NN_0}: 
\|b\|_{2-m}^2:
=\sum_{n=0}^\infty |b_n|^2 
(n!)^{2-m}<\infty\right\}.$$
\end{lemma}

Therefore, we have the  Gelfand triple
\begin{align}
\label{eq:6Mar16b}
\bigcap_{m=1}^\infty 
\cF_m\subset \cF_1\subset 
\bigcup_{m=1}^\infty \cF_{2-m}.
\end{align}

The inclusion map from $\cF_{m}$ into $\cF_{m+1}$ is nuclear, and it follows that 
$\cap_{m=1}^\infty\cF_{2-m}$ is a Fr\'echet nuclear space, and in particular a perfect space in the terminology of Gelfand and Shilov; 
see \cite{GS2_english}.  The dual space $\cup_{m=1}^\infty
\cF_{2-m}$ has two different set of properties, topological and algebraic; the first follow from the theory of perfect spaces, and the structure algebra 
comes from the form of the weights. The fact that the product is jointly continuous
comes from the theory of reflexive Fr\'echet spaces; see 
\cite[IV.26, Theorem 2]{MR83k:46003}.\\

We begin with the topological properties. Although not metrizable, the space $\cup_{m=1}^\infty\mathcal F_{2-m}$ behaves well with respect 
to sequences and compactness:
\begin{enumerate}
\item A sequence converges in the strong (or weak) topology of the dual if and only if its elements are in 
one of the spaces $\cF_{2-m}$ and converges in the topology of the latter; see \cite[p. 56]{GS2_english}.

\item A subset of $\cup_{m=1}^\infty\cF_{2-m}$ is compact in the strong topology of the dual if and only if it is included in one the spaces
$\cF_{2-m}$ and compact in the topology of the latter; see \cite[p. 58]{GS2_english}.

\end{enumerate}

These properties allow us to reduce to the Hilbert space setting and sequences the study of continuous functions from a compact metric space into $\cup_{m=1}^\infty\cF_{2-m}$.
\\

The algebra structure is given by the convolution product (or Cauchy product) defined as follows:

\begin{equation}
\label{convol}
a*b:=\left(\sum_{k=0}^n 
a_kb_{n-k}\right)_{n\in\NN_0},
\end{equation}
where $a=(a_n)_{n\in\mathbb N_0}$ and $b=(b_n)_{n\in\mathbb N_0}$ belong to the dual.

\begin{proposition}
The space 
\begin{align*}
\cF^\prime:=
\bigcup_{m=1}^\infty \cF_{2-m}
=\left\{ b=(b_n)_{n\in\NN_0}: 
\exists m\ge1,\|b\|_{2-m}:=
\sum_{n=0}^\infty 
\frac{|b_n|^2}{(n!)^{m-2}}<\infty 
\right\},
\end{align*}
is a topological algebra; the convolution product is jointly continuous with respect to the two variables, and
satisfies 
\begin{equation}
\label{vage1234}
\| a*b \|_{2-p}\le A(q-p)
\|a\|_{2-q}\|b\|_{2-p},
\end{equation}
for every
$a\in \cF_{2-q}$ 
and 
$b\in \cF_{2-p}$, 
where 
$p,q\in\NN$ 
such that 
$q\ge p+1$. 
\end{proposition}

The weights $\alpha_n=n!$ satisfy
\[
\alpha_{m+n}=\sqrt{(m+n)!}
\ge\sqrt{m!n!}=\alpha_m\alpha_n
\]
for every $m,n\in\NN_0$ and 
$\sum_{n=0}^\infty 
(\alpha_n)^{-2}=\sum_{n=0}^\infty 
\frac{1}{n!}=e<\infty$. Using these properties of the weight the statements in the proposition follow then from \cite{MR3404695} or, in a maybe more explicit way, 
from \cite[Exercise 5.4.8 p. 260-261]{CAPB_2}, with
\[
A(q-p)
=\left(\sum_{n=0}^\infty 
\alpha_n^{2(p-q)}\right)^{1/2}
=\left(\sum_{n=0}^\infty \left(\frac{1}{n!}
\right)^{q-p}\right)^{1/2}<\infty
\]
for $q-p\ge 1$.\\

We note that \eqref{vage1234} is called V\"age inequality, and originates with the work of V\"age; see \cite{MR2387368,vage96}.\\

Consider now a $\cF_1$-valued function, say $f$, defined a compact set (for instance $[0,1]$). When viewing $f$ as $\cup_{m=1}^\infty\cF_{2-m}$-valued, one can define 
differentiability and compute explicitly the derivative, which will take values in 
one of the spaces $\mathcal F_{2-m}$ rather than in the Fock space itself.
Using V\"age inequality one can also consider stochastic type integrals of the form
\[
\int_0^1 f(t)*g(t)dt,
\]
where $f$ and $g$ are continuous from $[0,1]$ into $\mathcal F^\prime$ as Riemann integrals. The image of $[0,1]$ under the function $f* g$ is then compact and
the integral is computed in one of the spaces $\mathcal F_{2-m}$. See \cite{aal2,aal3,MR3231624,
MR3615375} for similar arguments and applications, the latter
in the setting of quaternionic stochastic processes. Finally, we refer to  \cite{ACSS} for the 
study of the quaternionic Fock space and to \cite{diki2} for some of 
its generalizations in the quaternionic setting.
\\\\
\textbf{Acknowledgment.} 
Daniel Alpay thanks the Foster G. and Mary McGaw Professorship in
Mathematical Sciences, which supported this research. The authors would like
to thank Prof. Karol Penson for his helpful remarks related to the kernels
(\ref{eq:22May18a}) and the Meijer functions, and Prof. Dmitrii Karp for
pointing out the references \cite{Karp2,Karp1}.
 
\bibliographystyle{plain}
\def\cprime{$'$} \def\cprime{$'$} \def\cprime{$'$}
  \def\lfhook#1{\setbox0=\hbox{#1}{\ooalign{\hidewidth
  \lower1.5ex\hbox{'}\hidewidth\crcr\unhbox0}}} \def\cprime{$'$}
  \def\cprime{$'$} \def\cprime{$'$} \def\cprime{$'$} \def\cprime{$'$}
  \def\cprime{$'$}

\end{document}